\documentclass[final,nomarks]{dmtcs-episciences}

\usepackage{hyperref}
\usepackage{amsmath,amsfonts,amssymb,amsthm}
\usepackage[T1]{fontenc}
\usepackage{xcolor}
\newtheorem{theorem}{Theorem}
\newtheorem{conjecture}{Conjecture}

\newtheorem{remark}{Remark}

\DeclareMathOperator{\power}{Pow}

\usepackage[round]{natbib}

\author{Sebastian Czerwiński }
\title{On harmonious coloring of hypergraphs}%
\affiliation{%
  Institute of Mathematics, 
  University of Zielona Góra, Poland}

\publicationdata{vol. 26:2}{2024}{6}{10.46298/dmtcs.11101}{2023-03-21; 2023-03-21; 2023-11-08}{2024-04-05}
\begin{document}

\maketitle

\begin{abstract}
A \emph{harmonious coloring} of a $k$-uniform hypergraph $H$ is a vertex coloring such that no two vertices in the same edge have the same color,  and each $k$-element subset of colors appears on at most one edge. The \emph{harmonious number} $h(H)$ is the least number of colors needed for such a coloring. 
			
The paper contains a new proof of the upper bound $h(H)=O(\sqrt[k]{k!m})$ on the harmonious number of hypergraphs of maximum degree $\Delta$ with $m$ edges. We use the local cut lemma of A. Bernshteyn.
\end{abstract}
	
\section{Introducion}

	Let $H=(V,\mathcal{E})$ be a $k$-uniform hypergraph with the set of vertices $V$ and the set of edges $\mathcal{E}$. The set of edges is a family of $k$-element sets of $V$, where $k\geq 2$.
	
	A \emph{rainbow} coloring $c$ of a hypergraph $H$ is a map $c: V \mapsto\{1,\ldots,r\}$ in which no two vertices in the same edge have the same color. If two vertices are in the same edge $e$ with the same color, we say that the edge $e$ is \emph{bad}.
	
	A coloring $c$ is called \emph{harmonious} if $c(e)\neq c(f)$ for every pair of distinct edges $e,f\in\mathcal{E}$ and $c$ is a rainbow coloring.
	
	We say that distinct edges $e$ and $f$ have the \emph{same pattern of colors} if $c(e \setminus f)=c(f \setminus e)$ and there is no uncolored vertex in the set $e \setminus f$.
	
	Let $h(H)$ be the least number of colors needed for a harmonious coloring of $H$.
	In \cite{Bosek}, the following result is proved
	\begin{theorem}[\cite{Bosek}]
		For every $\varepsilon>0$ and every $\Delta>0$ there exist integers $k_0$ and $m_0$ such that every $k$-uniform hypergraph $H$ with $m$ edges (where $m\geq m_0$ and $k\geq k_0$) and maximum degree $\Delta$ satisfies
		\[h(H)\leq(1+\varepsilon)\frac{k}{k-1}\sqrt[k]{\Delta(k-1)k!m}.\]
	\end{theorem}
 
        \begin{remark}
        The paper \cite{Bosek} contains the following upper bound on the harmonious number 
        \[h(H)\leq\frac{k}{k-1}\sqrt[k]{\Delta(k-1)k!m}+1+\Delta^2+(k-1)\Delta+\sum_{i=2}^{k-1}\frac{i}{i-1}\sqrt[i]{(i-1)i\frac{(k-1)\Delta^2}{k-i}}.\]
        \end{remark}
        
	The proof of this theorem is based on the entropy compression method, see \cite{GrytczukKM, EsperetP}.
	
	Because a number $r$ of used colors must satisfy the inequality $\binom{r}{k}\leq m$, we get the lower bound $\Omega(\sqrt[k]{k!m})$. By these observations, it is conjectured by \cite{Bosek} that 
	\begin{conjecture}
		For each $k,\Delta\geq 2$ there exist a constant $c=c(k,\Delta)$ such  that every $k$-uniform hypergraph $H$ with $m$ edges and maximum degree $\Delta$ satisfies \[h(H)\leq \sqrt[k]{k!m}+c.\]
	\end{conjecture}
	A similar conjecture was posed by \cite{EdwardsDegree} for simple graphs $G$. He prove there that \[h(G)\leq(1+o(1))\sqrt{2m}.\]
	
	There are many results about the harmonious number of particular classes of graphs, see \cite{Aflaki, Akbari, Edwards1, Edwards2, Edwards0, Krasikov} or \cite{Aigner, Balister1, Balister2, Bazgan, Burris}.
 
	The present paper a new contains proof of the theorem of Bosek et al. We use a different method, the local cut lemma of \cite{Bernshteyn, Bernshteyn2}. The proof is simpler and shorter than the original proof of Bosek et al.
	
	\section{A special version of the local cut lemma}
	Let $A$ be a family of subsets of a powerset $\power(I)$, where $I$ is a finite set. We say that it is \emph{downwards-closed} if for each $S\in A$, $\power(S)\subseteq(A)$. A subset $\partial A$ of $I$ is called \emph{boundary} of a downwards-closed family $A$ if 
	\[\partial A:=\{i\in I:\ S\in A \text{ and } S\cup \{i\}\not\in A \text{ for some }S\subseteq I\setminus\{i\}\}.\]
	Let $\tau\colon T\mapsto [1;+\infty)$ be a function, then for every $X\subseteq I$ we denote by $\tau(X)$ a number
	\[\tau(X):=\prod_{x\in X} \tau(x).\]
	Let $B$ be a random event, $X\subseteq I$ and $i\in I$. We introduce two quantities:
	\[\sigma^{A}_{\tau}(B,X):=\max_{Z\subseteq I\setminus X}\Pr(B \text{ and } Z\cup X \not\in A|Z\in A)\cdot\tau(X)\]
	and
	\[\sigma^{A}_{\tau}(B,i):=\min_{i\in X\subseteq I}\sigma^{A}_{\tau}(B,X).\]
        If $\Pr(Z\in A)=0$, then $\Pr(P|Z\in A)=0$, for all events $P$.
	\begin{theorem}[\cite{Bernshteyn}]\label{th:lcl}
		Let $I$ be a finite set. Let $\Omega$ be a probability space and let $A\colon\Omega\mapsto\power(\power(I))$ be a random variable such that with probability 1, $A$ is a nonempty downwards-closed family of subsets of $I$. For each $i\in I$, Let $\mathcal{B}(i)$ be a finite collection of random events such that whenever $i\in\partial A$, at least one of the events in $\mathcal{B}(i)$ holds. Suppose that there is a function $\tau\colon I\mapsto [1,+\infty)$ such that for all $i\in I$ we have 
		\[\tau(i)\geq1+\sum_{B\in \mathcal{B}(i)} \sigma^{A}_{\tau}(B,i).\]
		Then $Pr(I\in A)\geq1/\tau(I)>0$.
	\end{theorem}
	
	\section{Proof of Theorem 1}
	We choose a coloring $f:V\mapsto \{1,\ldots,t\}$ uniformly at random. Let $A$ be a subset of the power set of $V$ given by
	\[A:=\{S\subseteq V:c \text{ is a harmonious coloring of }H(V)\}.\]
	It is a nonempty downwards-closed family with probability 1 (the empty set is an element of $A$)
	
	By a set $\partial A$, we denote the set of all vertices $v$ such that there is an element $X$ of $A$ such that the coloring $c$ is not a harmonious coloring of $X\cup\{v\}$. If the coloring $c$ is not a harmonious coloring, then there is a bad edge or there are two different edges with the same pattern of colors. So, we define for every $v \in V$, a collection $\mathcal{B}(v)$ as union of sets: 
	\[\mathcal{B}^1(v):=\{B_e:v\in e\in E(H) \text{ and } e \text{ is not proper colored}\},\]
	and for every $i\in\{0,1,\ldots,k-1\}$,
	\[\mathcal{B}^2_i(v):=\{B_{e,f}:v\in e,f\in E(H) \text{ and } c(e)=c(f),\ |e\setminus f|=i\}.\]
	That is $\mathcal{B}(v)=\mathcal{B}^1(v)\cup \bigcup_{i=1}^{k-1}\mathcal{B}^2_i(v)$.
	
	We assume that the event $B_e$ happens if and only if the edge $e$ is the bad edge and the event $B_{e,f}$ happens if and only if the edges $e$ and $f$ have the same pattern of colors.
	
	We also assume that a function $\tau$ is a constant function, that is, $\tau(v)=\tau\in[1,+\infty)$. This implies that for any subset $S$ of $V$, we have $\tau(S)=\tau^{|S|}$.
	
	Now, we must find an upper bound on 
	\[\sigma^A_{\tau}(B,v)=\min_{X\subseteq V:v\in X} \max_{Z\subseteq V\setminus X} \Pr(B \wedge Z \cup X \not\in A|Z\in A)\tau(X),\] 
	where $v\in V$ and $B\in\mathcal{B}(v)$. We will use an estimation 
    \[\sigma^A_{\tau}(B,v) \leq \max_{Z\subseteq V\setminus X} \Pr(B|Z\in A)\tau(X).\] Now, we consider two cases.\\
	
	\noindent Case 1: $B\in\mathcal{B}^1$, i.e. $B=B_e$.\\
	We choose as $X$ the set $\{e\}$. Because the colors of distinct vertices are independent, we get an upper bound $\sigma^A_{\tau}(B_e,v)\leq \Pr(B_e)\tau^k$ (events $B_e$ and $"Z\in A"$ are independent).
	The probability $\Pr(\overline{B_e})$, opposite to $\Pr(B_e)$, fulfills   \[\Pr(\overline{B_e})=1-\frac{t}{t}\cdot\frac{t-1}{t}\cdot\ldots\cdot\frac{t-k+1}{t}\geq 1-(1-\frac{k-1}{t})^{k-1}.\]
	By Bernoulli's inequality, we get 
	\[\Pr(\overline{B_e})\geq 1-(1-\frac{k-1}{t}\cdot (k-1))=\frac{(k-1)^2}{t}.\]
	So, $\Pr(B_e) \leq \frac{k^2}{t}$.\\
	
	\noindent Case 2: $B\in\mathcal{B}^2_i$, i.e. $B=B_{e,f}$ and $|e\setminus f|=i$.\\
	Now, we set $X=e \setminus f$. The probability $\Pr(B_{e,f})$ is bounded from above by $\frac{i!}{t^i}$. So, we get 
	\[\sigma^A_{\tau}(B_{e,f},v)\leq \Pr(B_{e,f})\tau^i \leq \frac{i!}{t^i}\tau^i.\]
	
	To end the proof we must find an upper bound on sizes of sets $\mathcal{B}^1(v)$, $\mathcal{B}^2_0(v)$ and $\mathcal{B}^2_i(v)$, where $ i>0$. Because the degree of a vertex is bounded by $\Delta$ and the number of edges is $m$ we get that 
	\[|\mathcal{B}^1(v)|\leq \Delta \text{ and } |\mathcal{B}^2_0(v)|\leq \Delta m.\]
	The hardest part is an upper bound on $\mathcal{B}^2_i(v),\ i>0$. 
	The number of edges $f$ such that $e\setminus f =i$ is bounded by $\frac{k\Delta}{k-i}$. There are at most $k\Delta$ edges with a nonempty intersection with the edge $e$, and the edge $f$ has exactly $k-i$ common elements with $e$. So, we have $|\mathcal{B}^2_i(v)|\leq \Delta \frac{k\Delta}{k-i}$.
	To apply Theorem \ref{th:lcl} we must find $\tau\in[1,+\infty)$ and $c\in\mathbb{N}$ such that, for all $v\in V$, the inequality below holds 
	\[\tau\geq 1+\Delta\frac{k^2}{t}\tau^k+\Delta m \frac{k!}{t^k}\tau^k+\sum_{i=1}^{k-1}\Delta\frac{k\Delta}{k-i}\frac{i!}{t^i}\tau^i.\]
	If we choose $\tau=\frac{k}{k-1}$ and $t=\frac{k}{k-1}\sqrt[k]{\Delta(k-1)k!m}(1+\varepsilon)$, it is easy to see that the inequality holds for sufficiently large $m$.
	
\section*{Acknowledgments}
I would like to thank the anonymous reviewer for his/her suggestions and comments.
\bibliographystyle{abbrvnat}
\bibliography{harm2}

\end{document}